\documentclass[11pt,a4paper]{article}
\usepackage[T1]{fontenc}
\usepackage{lmodern}
\usepackage[margin=30mm]{geometry}
\usepackage{amsmath,amssymb,amsthm}
\usepackage[hidelinks]{hyperref}
\hypersetup{pdftitle={A Counterexample to Single-Control Null Controllability of Stochastic Heat Equations},pdfauthor={Feizheng Huang and Zhaobo Liu}}
\newtheorem{theorem}{Theorem}
\newtheorem{corollary}{Corollary}
\newcommand{\E}{\mathbb E}
\newcommand{\Prob}{\mathbb P}
\newcommand{\F}{\mathcal F}
\newcommand{\R}{\mathbb R}
\newcommand{\T}{\mathbb T}
\newcommand{\ind}{\mathbf 1}
\newcommand{\norm}[1]{\lVert #1\rVert}
\newcommand{\ip}[2]{\langle #1,#2\rangle}
\newcommand{\dd}{\,\mathrm d}
\DeclareMathOperator{\dist}{dist}
\title{A Counterexample to Single-Control Null Controllability of Stochastic Heat Equations}
\author{%
\begin{tabular}{cc}
Feizheng Huang & Zhaobo Liu \\
{\small\href{mailto:2024640043@mails.szu.edu.cn}{\texttt{2024640043@mails.szu.edu.cn}}} &
{\small\href{mailto:liuzhaobo@szu.edu.cn}{\texttt{liuzhaobo@szu.edu.cn}}}
\end{tabular}\\[0.4em]
{\small Institute for Advanced Study, Shenzhen University, Shenzhen, China}
}
\date{}
\begin{document}
\maketitle
\begin{abstract}
We disprove the conjecture that a single localized drift control suffices for null controllability of stochastic parabolic equations with bounded adapted coefficients. On a bounded interval with homogeneous Dirichlet boundary conditions, we construct a deterministic, nonnegative smooth noise coefficient supported at positive distance from the control interval. One fixed nonzero smooth initial state cannot be driven to zero almost surely at any positive time by any progressively measurable control with finite expected $L^2$ energy. The proof uses products of Brownian increments in the adjoint terminal data to keep the initial pairing bounded away from zero while the observation norm tends to zero. On the flat torus, the construction yields counterexamples in every Gevrey class of order $s>1$ whenever the complement of the control set contains an open ball. Thus Gevrey regularity above order one cannot replace analyticity as a general sufficient condition for null controllability.
\end{abstract}

\section{Introduction}

For the heat equation on a smooth bounded connected domain, a control acting on any nonempty open subset can drive every initial state to zero in any positive time. L\"u and Zhang \cite[\S12.2]{LZ} conjectured the same property for stochastic parabolic equations with bounded adapted coefficients and a single localized drift control. S\^{\i}rbu and Tessitore \cite[\S5]{ST} had already identified spatially dependent noise as an open case. We give a negative answer for a deterministic smooth coefficient in one space dimension.

The equation in \cite[(12.2)]{LZ} is
\begin{equation}
\begin{aligned}
 \mathrm dy&=(\Delta y+a_1y+Bu)\dd t+a_2y\dd W_t
       &&\text{in }(0,T)\times G,\\
 y&=0 &&\text{on }(0,T)\times\partial G,\\
 y(0)&=y_0 &&\text{in }G .
\end{aligned}
\label{forward}
\end{equation}
Here $a_1,a_2$ are bounded, progressively measurable real coefficient fields, $G$ is a smooth bounded connected domain, and $G_0\subset G$ is nonempty and open. The operator $B:L^2(G_0)\to L^2(G)$ extends functions by zero. We take $W$ to be the coordinate Brownian motion on canonical Wiener space, with Wiener measure $\Prob$ and its usual augmented natural filtration $(\F_t)_{t\ge0}$. Let $U_T$ be the space of progressively measurable $L^2(G_0)$-valued controls satisfying
\[
 \norm{u}_{U_T}^2=\E\int_0^T\norm{u(t)}_{L^2(G_0)}^2\dd t<\infty.
\]
Null controllability at $T$ means that every deterministic $y_0\in L^2(G)$ admits $u\in U_T$ with $y(T)=0$ almost surely.

Tang and Zhang \cite{TZ} obtained null controllability with a second control in the diffusion. For a single control in the drift, L\"u \cite{Lu} proved null controllability for a class of equations with a spatially independent noise coefficient. Liu, Wu, Yang and Zhong \cite[\S3.3]{LWYZ} obtained observability from measurable sets for a spatially independent adapted noise coefficient. Other positive results assume that the noise operators are diagonal in the Dirichlet eigenbasis with a uniform bound on the noise intensity of the modes \cite[Theorem~1.10]{HLP}, or that the noise multiplier on the flat torus admits a bounded holomorphic extension with a common strip width and a time-integrated norm bound uniform over sample paths \cite[Theorem~1.5]{LW}.

For two coupled stochastic parabolic equations, X.~Liu \cite{Liu} exhibited loss of controllability under a change in the diffusion coupling. Our example uses a scalar equation. Its noise multiplier vanishes near the control set and equals one on another interval, so multiplication takes the positive first Dirichlet eigenfunction outside its eigenspace.

\section{The counterexample}

Fix the intervals
\[
 G=(0,10),\qquad G_0=(1,2),\qquad D=(6,7),\qquad
 D'=(25/4,27/4).
\]
All functions with compact support in these subintervals are extended by zero to $G$.

\begin{theorem}\label{thm:main}
Let $b\in C_c^\infty(D)$ satisfy $0\le b\le1$ and $b=1$ on $D'$. Choose any nonzero $g\in C_c^\infty(D')$ with $g\ge0$. In \eqref{forward} take $a_1=0$, $a_2=b$ and $y_0=g$. Then for every $T>0$ and every $u\in U_T$ the variational solution satisfies
\[
  \Prob\bigl(\norm{y(T)}_{L^2(G)}>0\bigr)>0 .
\]
\end{theorem}

\section{Proof of Theorem \ref{thm:main}}

Let $S(t)$ be the Dirichlet heat semigroup on the real space $H=L^2(G)$. The adjoint $B^*$ is restriction to $G_0$. On the torus use the periodic semigroup with coordinate balls $D'\Subset D$ at positive distance from $G_0$. Define the self-adjoint contractions
\[
 K(r)=\int_0^1 S\big((1-s)r\big)bS(sr)\dd s,\qquad r\ge0,
\]
where $b$ acts by multiplication. Write $\E_t=\E[\cdot\mid\F_t]$. For an uncontrolled state $Y$ started at time $t$, the mild formula gives $\E_tY(t+r)=S(r)Y(t)$, and conditional It\^o isometry gives $\E_t[(W_{t+r}-W_t)Y(t+r)]=rK(r)Y(t)$.

Take $0=t_0<\cdots<t_N=T$ with nonincreasing lengths $h_j=t_j-t_{j-1}\le1$, $1\le j\le N$. Set $h_{N+1}=0$, $q_{N+1}=g$, $q_{j-1}=K(h_j)q_j$ for $j=N+1,\ldots,1$, so $q_N=g$. Put $P_0=1$ and $P_j=\prod_{i=1}^j[(W_{t_i}-W_{t_{i-1}})/h_i]$ for $1\le j\le N$.

Since $\E\norm{P_Ng}_H^2=\norm{g}_H^2\prod_{i=1}^N h_i^{-1}<\infty$, the adjoint $\mathrm dz_N=-(\Delta z_N+bZ_N)\dd t+Z_N\dd W_t$ with terminal datum $P_Ng$ has a solution $(z_N,Z_N)$ \cite[Theorem~1.21]{LZ}. Both processes are progressively measurable, $z_N$ has continuous paths, and
\[
 \E\left[\sup_{0\le t\le T}\norm{z_N(t)}_H^2
       +\int_0^T\norm{Z_N(t)}_H^2\dd t\right]<\infty.
\]
Testing the duality identity \cite[proof of Theorem~4.2]{ST} against $\F_t$-measurable indicators gives its conditional form. For $t_{j-1}\le t\le t_j$ and $r=t_j-t$, restart $Y$ from any deterministic $Y(t)=f\in H$. Conditioning at the nodes gives $\ip{z_N(t)}{f}_H=\E_t[P_j\ip{q_j}{Y(t_j)}_H]$, hence
\begin{equation}
 z_N(t)=\frac{P_{j-1}}{h_j}
   \big[(W_t-W_{t_{j-1}})S(r)+rK(r)\big]q_j .
 \label{adjoint-form}
\end{equation}
Let $S_{D'}(h)$ denote the Dirichlet heat semigroup on $D'$, acting after restriction to $D'$ and followed by zero extension. Domain monotonicity and $b=1$ on $D'$ give $K(h)f\ge S_{D'}(h)f$ for $h\ge0$ and nonnegative $f\in H$. Iteration gives
\begin{equation}
 \ip{g}{z_N(0)}_H=\ip{g}{q_0}_H\ge c_T>0,
 \qquad c_T\triangleq\norm{S_{D'}(T/2)g}_{L^2(D')}^2.
 \label{mass}
\end{equation}

Set $\kappa=\dist(G_0,D)^2/4>0$, using Euclidean distance on the interval and the flat geodesic distance on the torus. Let $\ind_D$ denote multiplication by the indicator of $D$. The Davies--Gaffney estimate \cite[Theorem~3.3]{CS} for the Friedrichs extension of $-\Delta$ gives
$\norm{B^*S(s)\ind_D}_{H\to L^2(G_0)}\le e^{-\kappa/s}$ for $s>0$ in both geometries. In $S(r)K(h_{j+1})$ and $K(r)$, heat evolution after multiplication by $b$ lasts at most $r+h_{j+1}\le2h_j$. Since $\norm{q_i}_H\le\norm{g}_H$, both $\norm{B^*S(r)q_j}_{L^2(G_0)}$ and $\norm{B^*K(r)q_j}_{L^2(G_0)}$ are at most $\norm{g}_H e^{-\kappa/(2h_j)}$.
The increment $W_t-W_{t_{j-1}}$ is centered and independent of $P_{j-1}$. Using its variance $t-t_{j-1}$ and $\E P_{j-1}^2=\prod_{i=1}^{j-1}h_i^{-1}$ in \eqref{adjoint-form} gives
\begin{equation}
 \E\int_{t_{j-1}}^{t_j}\norm{B^*z_N(t)}_{L^2(G_0)}^2\dd t
 \le \norm{g}_H^2\Big(\frac12+\frac{h_j}{3}\Big)
       e^{-\kappa/h_j}\prod_{i=1}^{j-1}h_i^{-1}.
 \label{obs-step}
\end{equation}

For $\lambda_1\ge1$ define $\kappa(\lambda_{j+1}-\lambda_j)=\log(2\lambda_j)$ for $j\ge1$, so $\lambda_j\uparrow\infty$. Integral comparison and the recursion give
\[
 \sum_{i=1}^{j-1}\lambda_i^{-1}\ge \kappa\int_{\lambda_1}^{\lambda_j}\frac{\dd x}{x\log(2x)}\to\infty,
 \qquad
 e^{-\kappa\lambda_j}\prod_{i=1}^{j-1}\lambda_i=2^{1-j}e^{-\kappa\lambda_1}.
\]
Use lengths $h_j=\lambda_j^{-1}$, shortening the last interval to reach $T$. For fixed $T$, the partition, $N=N(\lambda_1)$, and the adjoint pair depend on $\lambda_1$. Shortening only decreases the last interval's bound in \eqref{obs-step}. Summing gives $\norm{B^*z_N}_{U_T}^2\le 2\norm{g}_H^2e^{-\kappa\lambda_1}\to0$ as $\lambda_1\to\infty$. If $u$ drove $g$ to zero, the terminal duality identity \cite[proof of Theorem~4.2]{ST} and \eqref{mass} would give $c_T\le\norm{u}_{U_T}\norm{B^*z_N}_{U_T}$, a contradiction. The same sequence violates the state observability estimate \cite[(12.3)]{LZ}.

\section{A consequence for Gevrey coefficients}

L\"u and Wang \cite[\S6]{LW} asked whether their analyticity assumption on the noise could be relaxed to Gevrey regularity of order greater than one. Corollary~\ref{cor:gevrey} gives counterexamples for every such order. On $\T^d=(\R/2\pi\mathbb Z)^d$, with integer $d\ge1$, use $B$ and $U_T$ as above with ambient domain $\T^d$. For $s\ge1$, let $G^s(\T^d)$ consist of the real smooth functions $b$ for which there are $C,\rho>0$ such that
\[
 \norm{\partial^\alpha b}_\infty
 \le C \rho^{|\alpha|}(|\alpha|!)^s
 \qquad\text{for every multi-index }\alpha .
\]
Here $|\alpha|=\alpha_1+\cdots+\alpha_d$.

\begin{corollary}\label{cor:gevrey}
Fix $d\ge1$ and a nonempty open $G_0\subset\T^d$ whose complement contains an open ball. In the periodic equation
$\mathrm dy=(\Delta y+Bu)\dd t+by\dd W_t$, every deterministic time-independent $b\in G^1(\T^d)$ gives null controllability at every $T>0$.

For each $s>1$ there exist a deterministic time-independent coefficient $b_s\in G^s(\T^d)$ with $0\le b_s\le1$ and a nonzero smooth initial state $g_s$ such that, for every $T>0$ and every $u\in U_T$, the solution with coefficient $b_s$ and initial state $g_s$ satisfies \mbox{$\Prob(\norm{y(T)}_{L^2(\T^d)}>0)>0$}.
\end{corollary}

\begin{proof}
The $G^1$ bound gives a bounded holomorphic extension to a sufficiently narrow periodic strip. Multiplying its supremum norm by $\sqrt{T}$ gives the bound in \cite[(5)]{LW}, so \cite[Theorem~1.5]{LW} gives a control in $U_T$.

For $s>1$, choose coordinate balls $D'\Subset D$ with $\dist(D,G_0)>0$. The Gevrey cutoff theorem \cite[Proposition~3.11 and Corollary~3.12]{Rainer} gives $0\le b_s\le1$ compactly supported in $D$ and equal to one on $D'$. Extend $b_s$ by zero to $\T^d$ and fix a nonzero nonnegative $g_s\in C_c^\infty(D')$. The periodic case of the main proof applies.
\end{proof}

\enlargethispage{\baselineskip}
\paragraph{AI disclosure.} OpenAI Codex was used to help develop the proof, search the literature, and prepare the manuscript. The authors take responsibility for its content. The Lean~4 formalization is at {\urlstyle{same}\url{https://github.com/guohaoji21/stochastic-heat-lean/tree/eb36cac854a3160f1edba61af74fdacf0c8e12ef}}.

\end{document}